\documentclass[twoside,11pt]{amsart}
\usepackage{amsmath,amssymb,amsbsy,amstext,amsxtra,latexsym}
\usepackage{graphics, epsfig}
\usepackage{hyperref}
\newtheorem{theorem}{Theorem}[section]
\newtheorem{lemma}{Lemma}[section]
\newtheorem{proposition}{Proposition}[section]

\def\mP{{\mathbb P}}
\def\mQ{{\mathbb Q}}
\def\mC{{\mathbb C}}

\def\mZ{{\mathbb Z}}
\def\cO{{\mathcal O}}
\def\cX{{\mathcal X}}

\def\cV{{\mathcal V}}
\def\n{\noindent}
\def\b{\bigskip}
\def\m{\medskip}
\def\t{\tilde}
\def\T{{\mathbb T}}
\def\Ext{\operatorname{Ext}}
\def\ker{\operatorname{ker}}

\def\log{\operatorname{log}}

\begin{document}

\title{A simply connected surface of general type \\
       with $p_g=0$ and $K^2=2$}

\author{Yongnam Lee and Jongil Park}

\address{Department of Mathematics, Sogang University,
         Sinsu-dong, Mapo-gu, Seoul 121-742, Korea}

\email{ynlee@sogang.ac.kr}

\address{Department of Mathematical Sciences, Seoul
         National University,
         San 56-1, Sillim-dong, Gwanak-gu, Seoul 151-747, Korea}

\email{jipark@math.snu.ac.kr}

\date{August 30, 2006, revised at July 24, 2007}



\begin{abstract}
 In this paper we construct a simply connected, minimal, complex surface
 of general type with $p_g=0$ and $K^2 =2$ using a rational blow-down surgery
 and a $\mQ$-Gorenstein smoothing theory.
\end{abstract}

\maketitle

\section{Introduction}
\label{sec-1}

\markboth{YONGNAM LEE AND JONGIL PARK}{A SIMPLY CONNECTED SURFACE
          OF GENERAL TYPE WITH $p_g=0$ AND $K^2=2$}

 One of the fundamental problems in the classification of complex surfaces
 is to find a new family of simply connected surfaces of general type with
 $p_g=0$.
 Surfaces with $p_g=0$ are interesting in view of Castelnuovo's criterion:
 An irrational surface with $q=0$ must have $P_2\ge 1$. This class
 of surfaces has been studied extensively by algebraic geometers
 and topologists for a long time. The details of history and examples
 of surfaces with $p_g=0$ are given in \cite{Do}.
 In particular,
 simply connected surfaces of general type with $p_g=0$ are little known
 and the classification is still open. Therefore
 it is a very important problem to find a new family of simply connected
 surfaces with $p_g=0$.
 Although a large number of non-simply connected complex surfaces of
 general type with $p_g =0$ have been known
 (\cite{BHPV}, Chapter VII),
 until now the only previously known simply connected,
 minimal, complex surface of general type with $p_g=0$
 was Barlow surface \cite{b}. Barlow surface has $K^2=1$. The
 natural question arises if there is a simply connected surface of
 general type with $p_g=0$ and $K^2\ge 2$.

 Recently, the second author constructed a new simply connected
 symplectic $4$-manifold with $b_2^+ =1$ and $K^2 =2$ using a
 rational blow-down surgery \cite{p2}. After this construction,
 it has been a very intriguing question whether such a symplectic
 $4$-manifold admits a complex structure. On the other hand,
 the compactification theory of a moduli space of surfaces of general
 type was established during the last twenty years. It was originally
 suggested in \cite{KSB} and it was established by Alexeev's proof
 for boundness~\cite{Al} and by the Mori program for threefolds
 (\cite{KM} for details).
 It is natural to expect the existence of a surface with special
 quotient singularities in the boundary of a compact moduli space.
 Hence it would be possible to construct a new interesting surface
 of general type  by using a singular surface
 with special quotient singularities.
 This was also more or less suggested by Koll\'ar
 in a recent paper~\cite{kol06}.

 The aim of this paper is to give an affirmative answer for
 the question above. Precisely, we construct a simply connected, minimal,
 complex surface of general type with $p_g=0$ and $K^2 =2$ by modifying
 Park's symplectic $4$-manifold in~\cite{p2}.
 Our approach is very different from other classical constructions such as a
 finite group quotient and a double covering,
 due to Godeaux, Campedelli, Burniat and others.
 These classical constructions are well
 explained in \cite{Re}. Our main techniques are rational blow-down surgery
 and $\mQ$-Gorenstein smoothing theory.
 These theories are briefly reviewed and developed in Section 2, 3 and 4.
 We first consider a special cubic pencil in $\mP^2$ and blow up many times
 to get a projective surface $\t Z$ which contains several disjoint chains
 of curves representing the resolution graphs of special quotient singularities.
 Then we contract these chains of curves from the surface $\t Z$ to produce
 a projective surface $X$ with five special quotient singularities.
 The details of the construction are given in Section 3.
 Using the methods developed in Section 2,
 we prove in Section 4 that the singular surface $X$ has a
 $\mQ$-Gorenstein smoothing. And we also prove in Section 5 that
 the general fiber $X_t$ of the $\mQ$-Gorenstein smoothing is a minimal
 surface of general type with $p_g=0$ and $K^2 =2$.
 Finally, applying the standard arguments about Milnor fibers
 (\cite{LW} \S5, or \cite{Man01} \S1),
 we prove that the general fiber $X_t$ is diffeomorphic to a simply
 connected symplectic $4$-manifold $\t Z_{15,9,5,3,2}$ which is obtained
 by rationally blowing down along five disjoint configurations
 in $\t Z$. The main result of this paper is the following

\begin{theorem}
\label{thm-main}
  There exists a simply connected, minimal, complex surface
  of general type with $p_g=0$ and $K^2 =2$.
\end{theorem}

 We provide more examples of simply connected surfaces of general type
 with $p_g=0$ and $K^2=2$ in Section 6.
 By using a different configuration, we are also able to construct
 a minimal, complex surface of general type with $p_g=0, K^2 =2$ and
 $H_1=\mZ_2$~\cite{LP}.
 Furthermore, by a small modification of the main construction
 and using the same techniques, we also construct simply connected,
 minimal, complex surfaces of general type with $p_g=0$ and $K^2 =1$.
 These constructions are explained in Appendix.

\m

\n{\em Acknowledgements}. The authors would like to thank Denis
 Auroux, Ronald Fintushel, Yoshio Fujimoto, Jonghae Keum,
 Keiji Oguiso and Noboru Nakayama for their valuable comments
 and suggestions to improve this article.
 The authors acknowledge Igor Dolgachev, J\'anos Koll\'ar and Miles
 Reid for their useful comments and information about related papers.
 The authors also thank Ki-Heon Yun for his drawing figures.
 Finally, the authors deeply appreciate Marco Manetti for his kind
 allowance to use his simple and direct proof of Theorem~\ref{thm-2.1}
 in Section 2. Yongnam Lee was supported by the grant (KRF-2005-070-C00005)
 from the KRF and he was supported by JSPS Invitation Fellowship
 Program while he was visiting RIMS at Kyoto University.
 He also thanks Shigefumi Mori for an invitation and for kind hospitality
 during his stay at RIMS. He also holds a joint appointment at
 CQUeST, Sogang University.
 Jongil Park was supported by grant No. R01-2005-000-10625-0 from
 the KOSEF and he also holds a joint appointment in the Research
 Institute of Mathematics, Seoul National University.

\b

\section{$\mQ$-Gorenstein smoothing}
\label{sec-2}

 In this section we develop a theory of $\mQ$-Gorenstein smoothing
 for projective surfaces with special quotient singularities,
 which is a key technical ingredient in our result.

\m

 \n {\bf Definition.} Let $X$ be a normal projective surface with quotient
 singularities. Let $\cX\to\Delta$ (or $\cX/\Delta$) be a flat family of
 projective surfaces over a small disk $\Delta$. The one-parameter
 family of surfaces $\cX\to\Delta$ is called a {\it
 $\mQ$-Gorenstein smoothing} of $X$ if it satisfies the following
 three conditions;

\n(i) the general fiber $X_t$ is a smooth projective surface,

\n(ii) the central fiber $X_0$ is $X$,

\n(iii) the canonical divisor $K_{\cX/\Delta}$ is $\mQ$-Cartier.

\m

A $\mQ$-Gorenstein smoothing for a germ of a quotient singularity
$(X_0, 0)$ is defined similarly. A quotient singularity which
admits a $\mQ$-Gorenstein smoothing is called a {\it singularity
of class T}.

\begin{proposition}[\cite{KSB, Man91, Wa}]
\label{pro-2.1}
 Let $(X_0, 0)$ be a germ of two dimensional quotient
 singularity. If $(X_0, 0)$ admits a $\mQ$-Gorenstein smoothing over
 the disk, then $(X_0, 0)$ is either a rational double point or a
 cyclic quotient singularity of type $\frac{1}{dn^2}(1, dna-1)$ for some
 integers $a, n, d$ with $a$ and $n$ relatively prime.
\end{proposition}

\begin{proposition}[\cite{KSB, Man91, Wa2}]
\label{pro-2.2}
\begin{enumerate}
 \item The singularities ${\overset{-4}{\circ}}$ and
 ${\overset{-3}{\circ}}-{\overset{-2}{\circ}}-{\overset{-2}{\circ}}-\cdots-
 {\overset{-2}{\circ}}-{\overset{-3}{\circ}}$ are of class $T$.
 \item If the singularity
 ${\overset{-b_1}{\circ}}-\cdots-{\overset{-b_r}{\circ}}$ is of
 class $T$, then so are
 $${\overset{-2}{\circ}}-{\overset{-b_1}{\circ}}-\cdots-{\overset{-b_{r-1}}
 {\circ}}- {\overset{-b_r-1}{\circ}} \quad\text{and}\quad
 {\overset{-b_1-1}{\circ}}-{\overset{-b_2}{\circ}}-\cdots-
 {\overset{-b_r}{\circ}}-{\overset{-2}{\circ}}.$$ \item Every
 singularity of class $T$ that is not a rational double point can be
 obtained by starting with one of the singularities described in
 $(1)$ and iterating the steps described in $(2)$.
\end{enumerate}
\end{proposition}

 Let $X$ be a normal projective surface with singularities of class
 $T$. Due to the result of Koll\'ar and Shepherd-Barron \cite{KSB},
 there is a $\mQ$-Gorenstein smoothing locally for each singularity
 of class $T$ on $X$ (see Proposition~\ref{pro-ksb}).
 The natural question arises whether this local $\mQ$-Gorenstein
 smoothing can be extended over the global surface $X$ or not.
 Roughly geometric interpretation is the following:
 Let $\cup V_{\alpha}$ be an open covering of $X$ such that
 each $V_{\alpha}$ has at most one singularity of class $T$.
 By the existence of a local $\mQ$-Gorenstein
 smoothing, there is a $\mQ$-Gorenstein smoothing
 $\cV_{\alpha}/\Delta$. The question is if these families glue to
 a global one. The answer can be obtained by figuring out the
 obstruction map of the sheaves of deformation
 $T^i_X=Ext^i_X(\Omega_X,\cO_X)$ for $i=0,1,2$.
 For example, if $X$ is a smooth surface,
 then $T^0_X$ is the usual holomorphic tangent sheaf $T_X$ and
 $T^1_X=T^2_X=0$. By applying the standard result of deformations
 \cite{LS, Pal} to a normal projective surface with quotient
 singularities, we get the following

\begin{proposition}[\cite{Wa}, \S 4]
\label{pro-2.3}
 Let $X$ be a normal
 projective surface with quotient singularities. Then
\begin{enumerate}
 \item The first order deformation space of $X$ is represented by
 the global Ext 1-group $\T^1_X=\Ext^1_X(\Omega_X, \cO_X)$. \item
 The obstruction lies in the global Ext 2-group
 $\T^2_X=\Ext^2_X(\Omega_X, \cO_X)$.
 \end{enumerate}
\end{proposition}

 Furthermore, by applying the general result of local-global spectral
 sequence of ext sheaves (\cite{Pal}, \S 3) to deformation theory
 of surfaces with quotient singularities so that
 $E_2^{p, q}=H^p(T^q_X) \Rightarrow \T^{p+q}_X$,
 and by $H^j(T^i_X)=0$ for $i, j\ge 1$, we also get

\begin{proposition}[\cite{Man91, Wa}]
\label{pro-2.4}
 Let $X$ be a normal projective surface with quotient singularities.
 Then
 \begin{enumerate}
 \item We have the exact sequence
 $$0\to H^1(T^0_X)\to \T^1_X\to \ker [H^0(T^1_X)\to H^2(T^0_X)]\to 0$$
 where $H^1(T^0_X)$ represents the first order deformations of $X$
 for which the singularities remain locally a product.
 \item If $H^2(T^0_X)=0$, every local deformation of
 the singularities may be globalized.
 \end{enumerate}
\end{proposition}

\n {\bf Remark}. The vanishing $H^2(T^0_X)=0$ can be obtained via
 the vanishing of $H^2(T_V(-\log \ E))$, where $V$ is the minimal resolution
 of $X$ and $E$ is the reduced exceptional divisors.

\m

 Let $V$ be a nonsingular surface and let $D$ be a simple
 normal crossing divisor in $V$. The short exact sequence
 $$0\to \cO_D(-D) \to \Omega_V|_D \to \Omega_D\to 0$$
 induces
 $$0\to T^0_D\to T_V|_D \to \cO_D(D) \to Ext^1(\Omega_D, \cO_D)\to 0 .$$
 Hence there is a short exact sequence
 $$0\to T^0_D\to T_V|_D \to \oplus_i N_{D_i|V}\to 0$$
 where $D_i$ is a component of $D$.
 Then the sheaf $T_V(-\log \ D)$ is
 defined by the following natural exact sequence
 $$0\to T_V(-\log \ D)\to T_V \to \oplus_i N_{D_i|V}\to 0 .$$

 \m

 M. Manetti provided the following lemma for us and he also suggested
 a simple and direct proof of Theorem~\ref{thm-2.1} below~\cite{Man}.
 Our original approach is through a theory of global smoothings
 of varieties with normal crossings developed by R. Friedman~\cite{Fr}.

\begin{lemma}[\cite{Man}]
\label{lem-2.1}
 Let $\pi: (V, E)\to (X_0, 0)$ be the minimal resolution of a germ
 of two dimensional quotient singularity and let
 $E=\cup E_i$ be the reduced exceptional divisor.
 Then  $R^1\pi_*T_V(-\log \ E)=R^2\pi_*T_V(-\log \ E)=0$.
\end{lemma}

\begin{proof}
 We may assume that $X_0$ is affine.
 For every effective divisor $Z$ supported in $E$, its tangent sheaf
 $\Theta_Z$ fits into the commutative diagram (\cite{BW}, p.70):
 \[ \begin{array}{ccccccccc}
 0 & \to & T_V(-\log \ E) & \to & T_V & \to & \oplus_i N_{E_i|V} & \to & 0\\
   &     & \downarrow    &  &   \downarrow & & \parallel & & \\
 0 & \to & \Theta_Z &   \to   & T_V\otimes\cO_Z & \to & \oplus_i N_{E_i|V} &
 \to & 0.
 \end{array}\]
 Snake lemma gives the exact sequence
 \[ 0 \to T_V(-Z) \to T_V(-\log \ E)\to \Theta_Z \to 0 \]
 and then, for $Z$ sufficiently big,
 we get $H^1(T_V(-\log \ E))=H^1(\Theta_Z)$: If $F$ is a locally
 free sheaf on $V$, then $H^2(V, F)=0$ and, since the singularity is
 rational, we also get $H^1(V, F(-Z))=0$.
 The proof of two facts follows easily from (\cite{BHPV}, p.93--95),
 and it is sufficient that, for every component $E_i$, the restriction
 of $F(-Z)$ is a direct sum $\oplus_i\cO_{E_i}(a_i)$ with all $a_i\ge 0$.
 According to the tautness of quotient singularities \cite{Lau},
 we have $H^1(\Theta_Z)=0$ for every $Z$ sufficiently large.
\end{proof}

\begin{theorem}
\label{thm-2.1}
 Let $X$ be a normal projective surface with
 singularities of class $T$.
 Let $\pi: V\to X$ be the minimal resolution and
 let $E$ be the reduced exceptional
 divisors. Suppose that $H^2(T_V(-\log \ E))=0$.
 Then there is a $\mQ$-Gorenstein smoothing of $X$.
\end{theorem}

\begin{proof}(\cite{Man})
 By Lemma~\ref{lem-2.1} above,
 $R^1\pi_*T_V(-\log \ E)=R^2\pi_*T_V(-\log \ E)=0$ and then,
 Leray spectral sequence implies
 \[ H^2(\pi_*T_V(-\log \ E))= H^2(T_V(-\log \ E))=0 .\]
 Since $\pi_*T_V=T^0_X$ \cite{BW}, we have
 \[ 0 \to \pi_*T_V(-\log \ E) \to T^0_X \to \Delta \to 0, \]
 where $\Delta$ is supported in the singular locus of $X$.
 Therefore we get
 \[ H^2(T_X^0)=H^2(\pi_*T_V(-\log \ E))=0.\]
 Hence every local deformation of the singularities may be globalized
 by Proposition~\ref{pro-2.4}.
 Furthermore, there is a $\mQ$-Gorenstein smoothing of $X$
 by Proposition~\ref{pro-ksb} below.
 \end{proof}

 Let $X$ be a normal projective surface with singularities of class
 $T$. Our concern is to understand $\mQ$-Gorenstein smoothings in
 $\T^1_X$, not the whole first order deformations. These special deformations
 can be constructed via local index one cover. Let $U\subset X$ be an
 analytic neighborhood with an index one cover $U'$. For the case of
 the field $\mC$, this index one cover is unique up to isomorphism.
 The first order deformations which associate $\mQ$-Gorenstein
 smoothings can be realized as the invariant part of $T^1_{U'}$.
 By the help of the birational geometry in threefolds and their
 applications to deformations of surface singularities,
 the following proposition is obtained.

\begin{proposition}[\cite{KSB}]
\label{pro-ksb}
 Let $(X_0, 0)$ be a germ of singularity of class T.
 Then $(X_0, 0)$ admits a $\mQ$-Gorenstein smoothing.
\end{proposition}

 Note that Theorem~\ref{thm-2.1} above can be
 easily generalized to any log resolution of $X$ by keeping the vanishing
 of cohomologies under blowing up at the points. It is obtained
 by the following well-known result.
 Proposition~\ref{pro-2.6} is also used in Section 4.

\m

\begin{proposition} [\cite{FZ}, \S1]
\label{pro-2.6}
 Let $V$ be a nonsingular surface and let $D$ be a
 simple normal crossing divisor in $V$.
 Let $f: V'\to V$ be a blowing up of $V$ at a point p of $D$.
 Set $D'=f^{-1}(D)_{red}$.
 Then $h^2(T_{V'}(-\log \ D'))=h^2(T_V(-\log \ D))$.
\end{proposition}
\m

\n{\bf Example}. We consider a pencil of cubics in $\mP^2$
 and blow up at the base points. Denote this surface by $Y$. Then
 $Y$ has an elliptic fibration over $\mP^1$. Assume that there is a
 nodal fiber. Let $\tau : Z\to Y$ be a blowing-up at the singular point
 on this nodal fiber. Let $F$ be the proper transform of the nodal fiber and
 let $E$ be the exceptional curve in $Z$.
 We construct a projective surface $X$ with one singularity of class
 $T$ by contracting -4-curve $F$.
 Since $H^2(Z, T_Z(-F))=0$, Theorem~\ref{thm-2.1} implies that
 $X$ has a $\mQ$-Gorenstein smoothing.
 Using the fact that $-K_X$ is an effective divisor together with
 the result of Manetti (\cite{Man91}, Theorem 21), one can also prove
 the existence of a $\mQ$-Gorenstein smoothing of $X$.
 Using the result again in~\cite{Man91},
 it is not hard to show that a general fiber $X_t$ of a $\mQ$-Gorenstein
 smoothing is a smooth rational elliptic surface with $K_{X_t}^2=0$.

\m

\n {\bf Remark}. Gompf \cite{Gom} constructed a symplectic
4-manifold by taking
 a fiber sum of two symplectic $4$-manifolds.
 To briefly recall Gompf's example, we start with a
 simply connected relatively minimal elliptic surface with a section
 and with $c_2=48$. There is only one up to diffeomorphism such an
 elliptic surface, which is called $E(4)$.
 It is also known that $E(4)$
 admits nine rational $(-4)$-curves as disjoint sections.
 Rationally blowing down $n$ $(-4)$-curves of $E(4)$ is the same as
 the normal connected sum of $E(4)$ with $n$ copies of $\mP^2$
 by identifying a conic in each $\mP^2$ with one $(-4)$-curve in $E(4)$.
 This $4$-manifold is denoted by $W_{4,n}$.
 The manifold $W_{4, 1}$ does not admit any complex structure because
 it violates the Noether inequality $p_g\le \frac{1}{2}K^2+2$ (cf. \cite{BHPV}).
 In fact, we have $H^2(E(4), T_{E(4)})\ne 0$. Therefore it does
 not satisfy the vanishing condition in Theorem~\ref{thm-2.1}:
 Let $h: E(4)\to\mP^1$ be an elliptic fibration.
 Assume $C$ is a general fiber of the map $f$.
 We have an injective map $0\to h^*\Omega_{\mP^1} \to \Omega_{E(4)}$
 and the map induces an injection
 $H^0(\mP^1, \Omega_{\mP^1}(2))\hookrightarrow H^0(E(4), \Omega_{E(4)}(2C))$
 by tensoring $2C$ on $0\to h^*\Omega_{\mP^1} \to \Omega_{E(4)}$.
 Since $K_{E(4)}=2C$,
 the cohomology $H^0(E(4), \Omega_{E(4)}(K_{E(4)}))$ is not zero.
 Hence the Serre duality implies that $H^2(E(4), T_{E(4)})$ is not zero.
 In fact, it is not hard to show that $h^2(E(4), T_{E(4)})=1$.

\b

\section{The main construction}
\label{sec-3}

 We begin with the rational elliptic surface
 $E(1)=\mP^2\sharp 9\overline{\mP}^2$.
 There are several ways to construct an elliptic fibration on $E(1)$.
 In this paper we use a special elliptic fibration $g: E(1) \to \mP^1$
 which is constructed as follows:
 Let $A$ be a line and $B$ be a smooth conic in $\mP^2$
 which are represented homologically by $h$ and $2h$ respectively,
 where $h$ denotes a generator of $H_2(\mP^2; \mZ)$.
 Choose another line $L$ in $\mP^2$ which meets $B$ at two distinct points
 $p, q$, and which also meets $A$ at a different point $r$.
 We may assume that the conic $B$ and the line $A$ meet at two different
 points which are not $p, q, r$.
 We now consider a cubic pencil in $\mP^2$ induced by $A+B$ and $3L$, i.e.
 $\lambda (A+B)+\mu (3L)$, for $[\lambda:\mu]\in \mP^1$.
 After we blow up first at three points $p, q, r$, blow up again three times
 at the intersection points of the proper transforms of $B, A$
 with the three exceptional curves $e_1, e_2, e_3$. Finally, blowing up again
 three times at the intersection points of the proper transforms of
 $B$ (representing homologically $2h-e_1-e_2-e_1'-e_2'$) and
 $A$ (representing homologically $h-e_3- e_3'$)
 with the three new exceptional curves $e_1', e_2', e_3'$,
 we get an elliptic fibration $E(1)=\mP^2\sharp 9\overline{\mP}^2$ over
 $\mP^1$. We denote the three new exceptional curves by $e_1'', e_2'', e_3''$,
 and let us denote this elliptic fibration by $g: Y=E(1) \to \mP^1$.
 We note that there is an $\tilde E_6$-singular fiber
 ($IV^*$ in Kodaira's table of singular elliptic fibers in~\cite{BHPV}, p.201)
 on the fibration $g: Y \to \mP^1$  which consists of the proper
 transforms of $L, e_1, e_1', e_2, e_2', e_3, e_3'$.
 We also note that there is one $I_2$-singular fiber
 (two rational -2-curves meeting two points) on $g: Y \to \mP^1$
 which consists of the proper transforms of the line $A$ and the conic $B$.
 Furthermore, by the proper choice of curves $A, B$
 and $L$ guarantees two more nodal singular fibers on $g: Y \to
 \mP^1$. For example, the pencil $\lambda(xy-z^2)z+\mu(x-y)^3$
 works. This pencil has singular fibers at
 $[\lambda: \mu] =[1:0], [0:1], [36:\sqrt3 i]$ and $[-36:\sqrt3 i]$.
 Hence the fibration $g: Y \to\mP^1$ has one $\tilde E_6$-singular
 fiber, one reducible $I_2$-singular fiber, and two nodal singular fibers
 (see Figure~\ref{fig-Y1} below).
 Notice that there are three sections $e_1'', e_2'', e_3''$ in $Y$,
 so that two sections $e_1'', e_2''$ meet the proper transform of the
 conic $B$ and the third section $e_3''$ meets the proper transform
 of the line $A$.

 \begin{figure}[hbtb]
 \begin{center}
 \setlength{\unitlength}{1mm}
 \includegraphics[height=3.5cm]{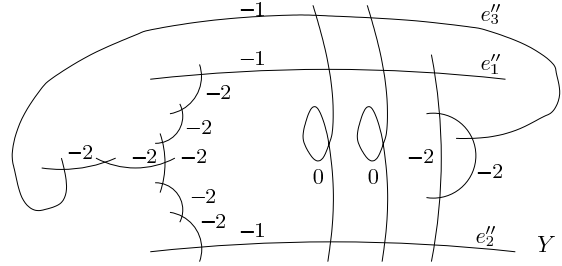}
 \end{center}
 \caption{A rational surface $Y$}
 \label{fig-Y1}
 \end{figure}

\m

\n {\bf Remark}. The existence of an elliptic fibration with
 fibres $I_1, I_1, I_2$ and $IV^*$ is known (\cite{Herf, Pers}).
 And, by the result in \cite{OSh}, this elliptic fibration has three
 sections which satisfy the configuration in Figure~\ref{fig-Y1}:
 Let $C_{1}$, $C_{2}$, $C_{3}$ be three ends of $IV^{*}$ and
 $D_{1}$, $D_{2}$ the two component of $I_{2}$.
 We number as the zero section $O$ passes $C_{1}$ and $D_{1}$.
 By the height formula, the generator $P$ of the Mordell Weil group passes
 $C_{2}$ (after changing $C_{2}$ and $C_{3}$ if necessary) and $D_{2}$.
 Then the section corresponding to $2P$ passes $C_{3}$ and $D_{1}$.
 Thus $O$, $2P$ and $P$ satisfy the configuration~\cite{Og}.

\m

\n{\bf Main Construction.} Let $Z:=Y\sharp 2{\overline \mP}^2$ be
 the surface obtained by blowing up at two singular points of two
 nodal fibers on $Y$, and denote this map by $\tau$. Then
 there are two fibers such that each consists of two $\mP^1$s,
 say $E_i$ and $F_i$, satisfying $E_i^2=-1$, $F_i^2=-4$ and
 $E_i\cdot F_i=2$ for $i=1,2.$
 Note that each $E_i$ is an exceptional curve and
 $F_i$ is the proper transform of a nodal fiber. Then the surface
 $Z$ has four special fibers; one $\tilde E_6$-singular
 fiber, one $I_2$-singular fiber consisting of $\tilde A$ and $\tilde B$
 which are the proper transforms of $A$ and of $B$, and
 two more singular fibers which are the union of $E_i$ and $F_i$ for $i=1, 2$.
 We denote the three sections $e''_1, e''_2, e''_3$ by $S_1, S_2, S_3$
 respectively. First, we blow up six times at the intersection points
 between two sections $S_1, S_2$ and $F_1, F_2, \tilde B$.
 It makes the self-intersection number of the proper transforms
 of $S_1$, $S_2$ and $\tilde B$ to be $-4$.
 We also blow up twice at the intersection points between the third section
 $S_3$ and $F_1, F_2$, so that the self-intersection number of the proper
 transform of $S_3$ is $-3$.

 Next, we blow up three times successively at the intersection point
 between the proper transform of $S_2$ and the exceptional curve in
 the total transform of $F_1$. It makes a chain of $\mP^1$,
 ${\overset{-7}{\circ}}-{\overset{-2}{\circ}}-{\overset{-2}{\circ}}
 -{\overset{-2}{\circ}}$,
 lying in the total transform of $F_1$.
 We also blow up three times successively at the intersection point
 between the proper transform of $S_2$ and the exceptional curve
 in the total transform of $F_2$,
 so that a chain of $\mP^1$,
 ${\overset{-7}{\circ}}-{\overset{-2}{\circ}}-{\overset{-2}{\circ}}
 -{\overset{-2}{\circ}}$,
 lies in the total transform of $F_2$. We note that the self-intersection
 numbers of the proper transforms of $F_1$ and $F_2$ are $-7$.
 Then we blow up at the intersection point
 between the proper transform of $S_1$ and the exceptional $-1$-curve
 intersecting the proper transform of $F_2$,
 so that it produces a chain of $\mP^1$,
 ${\overset{-2}{\circ}}-{\overset{-7}{\circ}}-{\overset{-2}{\circ}}
 -{\overset{-2}{\circ}} -{\overset{-2}{\circ}}$,
 lying in the total transform of $F_2$.
 Then we blow up again at the intersection point
 between the exceptional $-1$-curve and the rational $-2$-curve
 which is the right end of the above chain of $\mP^1$, so that
 it produces a chain of $\mP^1$,
 ${\overset{-2}{\circ}}-{\overset{-7}{\circ}}-{\overset{-2}{\circ}}-
 {\overset{-2}{\circ}}-{\overset{-3}{\circ}}-{\overset{-1}{\circ}}
 -{\overset{-2}{\circ}}$.
 We note that the self-intersection numbers of the proper transforms
 of $S_1$ and $S_2$ go to $-5$ and $-10$ respectively.

\begin{figure}[hbtb]
 \begin{center}
 \setlength{\unitlength}{1mm}
 \includegraphics[height=4cm]{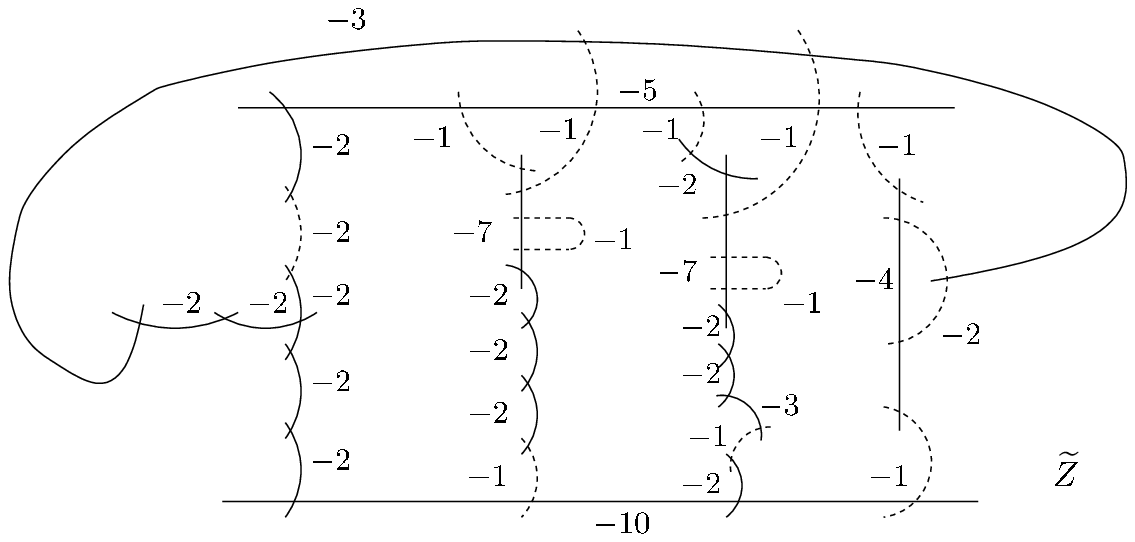}
 \end{center}
 \caption{A rational surface $\t Z$}
 \label{fig-tZ1}
 \end{figure}

 Next, we have a rational surface $\t Z:= Y \sharp 18{\overline \mP}^2$
 which contains five disjoint linear chains of $\mP^1$:
 ${\overset{-2}{\circ}}-{\overset{-10}{\circ}}
 -{\overset{-2}{\circ}}-{\overset{-2}{\circ}}-{\overset{-2}{\circ}}
 -{\overset{-2}{\circ}}-{\overset{-2}{\circ}}-{\overset{-3}{\circ}}$
 (which contains the proper transforms of two sections $S_2$, $S_3$
 and the five rational $-2$-curves in $\tilde E_6$-singular fiber),
 ${\overset{-7}{\circ}}-{\overset{-2}{\circ}}-{\overset{-2}{\circ}}
 -{\overset{-2}{\circ}}$,
 ${\overset{-2}{\circ}}-{\overset{-7}{\circ}}
 -{\overset{-2}{\circ}}-{\overset{-2}{\circ}}-{\overset{-3}{\circ}}$,
 ${\overset{-4}{\circ}}$
 and ${\overset{-5}{\circ}}-{\overset{-2}{\circ}}$
 (which contains the proper transforms of the section
 $S_1$ and the one rational $-2$-curve in $\tilde E_6$-singular fiber)
 (Figure~\ref{fig-tZ1}).

 Finally, we contract these five disjoint chains of $\mP^1$ from $\t Z$.
 Since it satisfies the Artin's criterion, it produces a projective
 surface with five singularities of class $T$ (\cite{Artin}, \S 2).
 We denote this surface by $X$. In Section 4 we will prove that
 $X$ has a $\mQ$-Gorenstein smoothing.
 And we will also show in Section 5 that a general fiber of
 the $\mQ$-Gorenstein smoothing is a simply connected, minimal,
 complex surface of general type with $p_g=0$ and $K^2=2$.

\m

 In the remaining of this section,
 we investigate a rational blow-down manifold of the surface $\t Z$.
 First we describe topological aspects of a rational blow-down surgery
 (\cite{fs, p1} for details):
 For any relatively prime integers $p$ and $q$ with $p > q > 0$,
 we define a configuration $C_{p,q}$ as a smooth $4$-manifold obtained
 by plumbing disk bundles over the $2$-sphere instructed by the following
 linear diagram
 $$\underset{u_{k}}{\overset{-b_k}{\circ}}-\underset{u_{k-1}}
 {\overset{-b_{k-1}}{\circ}}-\cdots -\underset{u_2}{\overset{-b_{2}}{\circ}}
 -\underset{u_1}{\overset{-b_{1}}{\circ}}$$

\n where $\frac{p^{2}}{pq-1} =[b_{k},b_{k-1}, \ldots, b_{1}]$ is the
 unique continued fraction with all $b_{i} \geq 2$, and  each
 vertex $u_{i}$ represents a disk bundle over the $2$-sphere whose
 Euler number is $-b_{i}$. Orient the $2$-spheres in $C_{p,q}$ so that
 $u_{i}\cdot u_{i+1} = +1$. Then the configuration $C_{p,q}$
 is a negative definite simply connected smooth $4$-manifold whose boundary
 is the lens space $L(p^2, 1-pq)$.

 \m

\n {\bf Definition.} Suppose $M$ is a smooth $4$-manifold
 containing a configuration $C_{p,q}$. Then we construct
 a new smooth $4$-manifold $M_p$,
 called a {\em $($generalized$)$ rational blow-down} of $M$,
 by replacing $C_{p,q}$ with the rational ball $B_{p,q}$.
 Note that this process is well-defined, that is, a new smooth
 $4$-manifold $M_{p}$ is uniquely determined
 (up to diffeomorphism) from $M$ because each diffeomorphism
 of $\partial B_{p,q}$ extends over the rational ball $B_{p,q}$.
 We call this a {\em rational blow-down} surgery.
 Furthermore, M. Symington proved that a rational blow-down manifold
 $M_{p}$ admits a symplectic structure in some cases.
 For example, if $M$ is a symplectic $4$-manifold containing
 a configuration $C_{p,q}$ such that all $2$-spheres $u_i$ in $C_{p,q}$
 are symplectically embedded and intersect positively,
 then the rational blow-down manifold $M_{p}$ also admits a symplectic
 structure~\cite{sy1, sy2}.

\m

 Now we perform a rational blow-down surgery on the surface $\t Z$
 constructed in the main construction.
 Note that the surface $\t Z$ contains five disjoint configurations
 - $C_{15,7},\, C_{9,4},\, C_{5,1},\, C_{3,1}$ and $C_{2,1}$.
 Let us decompose the surface $\t Z$ into
 \[{\t Z}={\t Z}_{0}\cup\{C_{15,7}\cup C_{9,4} \cup C_{5,1} \cup C_{3,1}
 \cup C_{2,1}\}.\]
 Then the $4$-manifold, say  ${\t Z}_{15,9,5,3,2}$, obtained by
 rationally blowing down along the five configurations can be decomposed into
 \[{\t Z}_{15,9,5,3,2}=
 {\t Z}_{0}\cup\{B_{15,7}\cup B_{9,4} \cup B_{5,1} \cup B_{3,1}
 \cup B_{2,1}\}, \]
 where $B_{15,7},\, B_{9,4},\, B_{5,1},\, B_{3,1}$ and $B_{2,1}$
 are the corresponding rational balls.  We claim that

\begin{theorem}
\label{thm-3.1}
  The rational blow-down ${\t Z}_{15,9,5,3,2}$ of the surface
  $\t Z$ is a simply connected closed symplectic $4$-manifold
  with $b_2^+=1$ and $K^2 =2$.
\end{theorem}

\begin{proof}
 Since all the curves lying in the configurations
 $C_{15,7},\, C_{9,4},\, C_{5,1},\, C_{3,1}$ and $C_{2,1}$
 are symplectically (in fact, holomorphically) embedded $2$-spheres,
 Symington's result~\cite{sy1, sy2} guarantees the existence of
 a symplectic structure on the rational blow-down $4$-manifold
 ${\t Z}_{15,9,5,3,2}$.
 Furthermore, it is easy to check that
 $b_2^+({\t Z}_{15,9,5,3,2})=b_2^+(\t Z)=1$ and
 $K^2({\t Z}_{15,9,5,3,2}) = K^2(\t Z) + 20 = 2$.

 It remains to prove the simple connectivity of ${\t Z}_{15,9,5,3,2}$:
 Since $\pi_{1}(\partial B_{p,q}) \rightarrow \pi_{1}(B_{p,q})$
 is surjective (\cite{LW}, \S5), by Van-Kampen theorem,
 it suffices to show that $\pi_{1}({\t Z}_0) =1$.
 First, note that $\t Z$ and all four configurations $C_{15,7}$, $C_{9,4}$,
  $C_{5,1}$, $C_{3,1}$ and $C_{2,1}$ are all simply connected.
 Hence, applying Van-Kampen theorem on $\t Z$ inductively, we get
 \begin{eqnarray}
 1 = \pi_{1}({\t Z}_0)/< N_{i_{*}(\alpha)}, N_{{j_1}_{*}(\beta_1)},
    N_{{j_2}_{*}(\beta_2)}, N_{{j_3}_{*}(\beta_3)}, N_{{k}_{*}(\gamma)}> .
 \end{eqnarray}
 Here $i_{*}$, ${j_{1}}_{*}$, ${j_{2}}_{*}$, ${j_{3}}_{*}$ and ${k}_{*}$
 are induced homomorphisms by inclusions
 $i: \partial C_{3,1} \rightarrow {\t Z}_0$,
 $j_1: \partial C_{5,1} \rightarrow {\t Z}_0$,
 $j_2: \partial C_{9,4} \rightarrow {\t Z}_0$,
 $j_3: \partial C_{2,1} \rightarrow {\t Z}_0$ and
 $k: \partial C_{15,7} \rightarrow {\t Z}_0$ respectively.
 We may also choose the generators, say $\alpha$, $\beta_1$,
 $\beta_2$, $\beta_3$ and $\gamma$, of
 $\pi_{1}(\partial C_{3,1}) \cong {\mZ}_{9}$,
 $\pi_{1}(\partial C_{5,1}) \cong {\mZ}_{25}$,
 $\pi_{1}(\partial C_{9,4}) \cong {\mZ}_{81}$,
 $\pi_{1}(\partial C_{2,1}) \cong {\mZ}_{4}$ and
 $\pi_{1}(\partial C_{15,7}) \cong {\mZ}_{225}$,
 so that $\alpha$, $\beta_1$, $\beta_2$, $\beta_3$ and $\gamma$
 are represented by circles $\partial C_{3,1} \cap E'_1$
 (equivalently $\partial C_{3,1} \cap E'_2$ or $\partial C_{3,1} \cap E'_3$),
 $\partial C_{5,1} \cap E'_1$,
 $\partial C_{9,4} \cap E'_2$ (equivalently $\partial C_{9,4} \cap E''_2$),
 $\partial C_{2,1} \cap E'_3$ and
 $\partial C_{15,7} \cap E''_2$, respectively,
 where $E'_1$, $E'_2$, $E'_3$ and $E''_2$ are exceptional curves connecting
 the last $2$-spheres in the configurations $C_{3,1}$ and $C_{5,1}$,
 $C_{3,1}$ and $C_{9,4}$, $C_{3,1}$ and $C_{2,1}$,
 $C_{15,7}$ and $C_{9,4}$, respectively.
 Note that the circle cut out by a 2-sphere which intersects transversely
 one of the two
 end $2$-spheres in the configurations $C_{p, q}$ is a generator
 of $\pi_1$ of the lens space, and other circles cut out by a $2$-sphere which
 intersects
 transversely one of the middle $2$-spheres in the configurations $C_{p, q}$,
 is a power of the generator \cite{Mu}.
 Finally $N_{i_{*}(\alpha)}$, $N_{{j_i}_{*}(\beta_i)}$ and
 $N_{k_{*}(\gamma)}$ denote the least normal subgroups
 of $\pi_{1}({\t Z}_0)$ containing
 $i_{*}(\alpha)$, ${j_i}_{*}(\beta_i)$ and ${k}_{*}(\gamma)$ respectively.
 Note that there is a relation between
 $i_{*}(\alpha)$ and ${j_i} _{*}(\beta_i)$ when we restrict them to
 ${\t Z}_{0}$. That is, they satisfy either
 $i_{*}(\alpha) = \tau^{-1} \cdot {j_i}_{*}(\beta_i)\cdot \tau$ or
 $i_{*}(\alpha) = \tau^{-1} \cdot {j_i}_{*}(\beta_i)^{-1} \cdot \tau$
 (depending on orientations) for some path $\tau$,
 because one is homotopic to the other in
 $E'_i \setminus \{\mathrm{two\, \, open\, \, disks} \} \subset {\t Z}_{0}$.
 Hence, by combining two facts above, we get
 ${i_{*}(\alpha)}^{9} = (\tau^{-1} \cdot {j_3}_{*}(\beta_3)^{\pm 1} \cdot
 \tau)^{9} = \tau^{-1} \cdot {j_3}_{*}(\beta_3)^{\pm 9} \cdot \tau
  = 1 = {{j_3}_{*}(\beta_3)}^4$. Since the two numbers $9$ and $4$ are
  relatively prime, the element ${{j_3}_{*}(\beta_3)}$ should
  be trivial. So the relation
  $i_{*}(\alpha) = \tau^{-1} \cdot {j_3}_{*}(\beta_3)^{\pm 1} \cdot \tau$
  implies the triviality of $i_{*}(\alpha)$.
  And the triviality of $i_{*}(\alpha)$ also implies that
  ${{j_1}_{*}(\beta_1)}$ and ${{j_2}_{*}(\beta_2)}$  are trivial.
  Furthermore, since ${j_2}_{*}(\beta_2)$ and
  $k_{*}(\gamma)$ are also conjugate each  other, $k_{*}(\gamma)$ is
  also trivial. Hence,
  all normal subgroups $N_{i_{*}(\alpha)}$, $N_{j_{i_{*}}(\beta_i)}$
  and $N_{k_{*}(\gamma)}$ are trivial, so that relation (1)
  implies $\pi_{1}({\t Z}_0) =1$.
\end{proof}

\m

\n {\bf Remark}. In fact, one can prove that the rational blow-down
 $4$-manifold ${\t Z}_{15,9,5,3,2}$ constructed in Theorem~\ref{thm-3.1}
 above is not diffeomorphic to a rational surface
 $\mP^2\sharp 7\overline{\mP}^2$ by using a technique in~\cite{p2}.
 Furthermore, one can also prove that the symplectic $4$-manifold
 ${\t Z}_{15,9,5,3,2}$ is minimal by using a similar technique
 in~\cite{os}.

\b

\section{Existence of smoothing}
\label{sec-4}

 In this section we prove the existence of a $\mQ$-Gorenstein
 smoothing for the main example constructed in Section 3.

\begin{lemma}
\label{lem-4.1}
 Let $Y$ be a rational elliptic surface. Let $C$ be a
 general fiber of the elliptic fibration $g: Y\to\mP^1$.
 Then the global sections $H^0(Y, \Omega_Y(kC))$ are coming from
 the global sections $H^0(Y, g^*\Omega_{\mP^1}(k))$.
 In particular, $h^0(Y, \Omega_Y(kC))=k-1$ for $k\ge 1$.
\end{lemma}

\begin{proof}
 We have an injective map $0\to g^*\Omega_{\mP^1} \to \Omega_Y,$
 and the map induces an injection $H^0(\mP^1,
 \Omega_{\mP^1}(k))\hookrightarrow H^0(Y, \Omega_Y(kC))$ by tensoring
 $kC$ on $0\to g^*\Omega_{\mP^1} \to \Omega_Y$.

 The dimension of the cohomology $H^0(Y, \Omega_Y(kC))$ is $k-1$:
 Consider the commutative diagram of standard exact sequences
 $$\begin{array}{ccc}
 & 0 & \\
 & \downarrow & \\
 & \cO_{C}((k-1)C) & \\
 & \downarrow & \\
 0\to\Omega_Y((k-1)C)\to\Omega_Y(kC)\to & \Omega_Y(kC)|_{C} & \to 0\\
 & \downarrow & \\
 & \Omega_{C}(kC) & \\
 & \downarrow & \\
 & 0 & .\end{array}$$ The map from the long exact sequence of
 cohomologies from the vertical sequence
 $$H^0(C, \Omega_C(kC))\simeq H^0(C, \Omega_C) \to H^1(C,
 \cO_{C}((k-1)C))\simeq H^1(C, N_{C|Y}^\vee)$$ is the dual of the
 Kodaira-Spencer map $H^0(C, N_{C|Y})\to H^1(C, T_{C})$ because
 $\Omega_{C}=\cO_{C}$. This map is not zero since the fibration is
 non-trivial. Therefore we have that $H^0(Y,
 \Omega_Y(kC)|_{C})=H^0(C, \cO_{C})$ and $h^0(Y, \Omega_Y(kC))\le
 h^0(Y, \Omega_Y((k-1)C))+1$.
 Furthermore, the map $H^0(C, \cO_{C})\to H^1(Y, \Omega_Y)$ is injective
 because it factors through the exact sequence
 $$0\to\Omega_Y\to\Omega_Y(\log\, C)\to\cO_{C}\to 0,$$
 and the map $H^0(C, \cO_{C})\to H^1(Y, \Omega_Y)$
 is the first Chern class map.
 Therefore we obtain the vanishing $H^0(Y, \Omega_Y(C))=0$
 and $h^0(Y, \Omega_Y(kC))\le k-1$.
 Since $H^0(Y, \Omega_Y(kC))$ contains $H^0(Y, g^*(\Omega_{\mP^1}(k)))$,
 the dimension of $H^0(Y, \Omega_Y(kC))$ is $k-1$. It implies that all
 global sections of $\Omega_Y(kC)$ are coming from the global sections
 of $g^*(\Omega_{\mP^1}(k))=g^*(\cO_{\mP^1}(k-2))$.
\end{proof}

\n {\bf Example.}
 Let $Y$ be a smooth rational elliptic surface
 constructed by blowing up the base points of a pencil of cubics.
 Assume that there are two nodal fibers.
 Let $\tau : Z\to Y$ be a blowing-up at the singular points $p, q$ on
 two nodal fibers.
 Let $F_1, F_2$ be the proper transforms of two nodal fibers and
 let $E_1, E_2$ be the two exceptional divisors in $Z$.
 Then we contract two -4-curves $F_1, F_2$.
 It produces a projective surface $X$ with two
 singularities of class $T$ whose resolution graphs are both
 ${\overset{-4}{\circ}}$. Let $f : Z\to X$ be the contraction morphism.
 Then $X$ has a $\mQ$-Gorenstein
 smoothing: By Theorem~\ref{thm-2.1},
 it is enough to prove $H^2(Z, T_Z(-F_1-F_2))=0$.
 It is equal to prove $H^0(Z, \Omega_Z(K_Z+F_1+F_2))=0$ by the Serre duality.
 Since $-\tau^*K_Y=F_i+2E_i$ for $i=1$ and $2$, we have
 $$2K_Z=\tau^*(2K_Y)+2E_1+2E_2=-F_1-F_2.$$
 Therefore
 $H^0(Z, \Omega_Z(K_Z+F_1+F_2))=H^0(Z, \Omega_Z(-K_Z))=
 H^0(Z, \Omega_Z(\tau^*C-E_1-E_2))$,
 where $C$ is a general fiber of an elliptic fibration form $Y$ to
 $\mP^1$. Since $H^0(Z, \Omega_Z(\tau^*C-E_1-E_2))\subseteq
 H^0(Z, \Omega_Z(\tau^*C))=H^0(Y, \Omega_Y(C))$, it is zero by
 Lemma 4.1. It is not hard to see that a general fiber $X_t$ of a
 $\mQ$-Gorenstein smoothing is a minimal Enriques surface.

\m

 If we contract more than two proper transforms of nodal fibers in a rational
 elliptic surface,
 then it does not satisfy the vanishing condition in Theorem~\ref{thm-2.1}.

\begin{proposition}
\label{pro-4.1}
 Let $Y$ be a smooth rational elliptic surface. Assume that the
 elliptic fibration $g: Y\to\mP^1$ is relatively minimal without
 multiple fibers and there are $j$ nodal fibers with $j\ge 3$.
 Let $\tau : Z\to Y$ be a blowing-up at the singular points $p_1,\ldots, p_j$
 on nodal fibers.
 Let $F_1, \ldots, F_j$ be the proper transforms of nodal fibers and
 let $E_1, \ldots, E_j$ be the exceptional divisors in $Z$.
 Then $H^2(Z, T_Z(-F_1-\cdots -F_j))\ne 0$.
\end{proposition}

\begin{proof}
 Let $h: Z\to \mP^1$ be an elliptic fibration, and let
 $s_1,\ldots, s_j$ be points on $\mP^1$ corresponding to nodal fibers
 $F_1, \ldots, F_j$. Then there is an injective map (\cite{EV}, p.80--81)
 $$0\to h^*\Omega_{\mP^1}(s_1+\cdots+ s_j)\to \Omega_Z(\log
 (F_1+\cdots + F_j+E_1+\cdots +E_j)).$$
 Let $C'$ be a general fiber of the map $h$, and $s$ be a point on $\mP^1$
 corresponding to a fiber $C'$. Since $K_Z=-C'+E_1+\cdots+E_j$,
 $h^0(Z, \Omega_Z(\log (F_1+\cdots + F_j+E_1+\cdots +E_j))(K_Z))\ge h^0(Z,
 h^*\Omega_{\mP^1}(s_1+\cdots+ s_j)(-s))=j-2.$
 Therefore, if $j\ge 3$, the Serre duality implies that
 $H^2(Z, T_Z(-\log (F_1+\cdots + F_j+E_1+\cdots +E_j)))\ne 0$.
 Since $E_i^2=-1$ for $i=1, \ldots, j$,
 it is also equal to $H^2(Z, T_Z(-\log (F_1+\cdots + F_j)))\ne 0$.
 Hence it implies $H^2(T_Z(-F_1-\cdots -F_j))\ne 0$.
\end{proof}

 By Proposition~\ref{pro-4.1} above, we need to choose a smooth rational
 elliptic surface with special fibers. This is one of the reasons
 why we consider a rational elliptic surface with four special singular
 fibers in Section 3.

\begin{lemma}
\label{lem-4.2}
 Let $Z=Y\sharp 2{\overline \mP}^2$ be the rational
 elliptic surface in the main
 construction. Let $F$ be the proper transform of the conic $B$ in $Z$,
 and $E$ be the proper transform of the line $A$ in $Z$. Let
 $D$ be the reduced subscheme of the $\tilde E_6$-singular
 fiber. Assume that $D$ is not whole $\tilde E_6$-singular
 fiber as a reduced scheme. Then $H^2(Z, T_Z(-F_1-F_2-F-D))=0$ and $H^2(Z,
 T_Z(-F_1-F_2-E-D))=0$.
\end{lemma}

 \begin{proof}
 By the Serre duality, it is equal to prove $H^0(Z,
 \Omega_Z(K_Z+F_1+F_2+F+D))=0$. Let $C$ be a general fiber in the elliptic
 fibration $g: Y\to\mP^1$. Since $K_Z=\tau^*(-C)+E_1+E_2$ and
 $\tau^*(C)=F_1+2E_1$, $H^0(Z,
 \Omega_Z(K_Z+F_1+F_2+F+D))\subseteq H^0(Z,
 \Omega_Z(\tau^*(C)+F+D))$. Since $F$ and $D$ are not changed by
 the map $\tau$, we have the same curves in $Y$. Then $H^0(Z,
 \Omega_Z(\tau^*(C)+F+D))=H^0(Y,
 \Omega_Y(C+F+D))$ by the projection formula.  We note that
 $\tau_*\Omega_Z=\Omega_Y$. Then the cohomology $H^0(Y,
 \Omega_Y(C+F+D))$ vanishes: We note that $H^0(Y,\Omega_Y(C+F+D))=H^0(Y,
 \Omega_Y(3C-E-D'))$ with $D'+D =\tilde E_6$-singular fiber.
 By Lemma~\ref{lem-4.1}, all global sections of
 $\Omega_Y(3C)$ are coming form the global sections of
 $g^*(\Omega_{\mP^1}(3))=g^*(\cO_{\mP^1}(1))$.
 But, if this global
 section vanishes on $E$ and $D'$ which lie on two different fibers,
 then it should be zero. Note that the dualizing sheaf of each fiber
 of the elliptic fibration is the structure sheaf of the fiber by
 using the adjunction formula. Therefore we have the vanishing
 $H^2(Z, T_Z(-F_1-F_2-F-D))=0$.
 The vanishing $H^2(Z, T_Z(-F_1-F_2-E-D))=0$ is obtained by the same proof.
 \end{proof}

\smallskip

\begin{theorem}
\label{thm-4.1}
 The projective surface $X$ with five singularities of class T
 in the main construction has a $\mQ$-Gorenstein smoothing.
\end{theorem}

\begin{proof}
 Let $D$ be the reduced scheme of the $\tilde E_6$-singular fiber
 minus the rational $-2$-curve $J$ in the main construction in Section 3.
 Note that the curve $J$ is not contracted from $\tilde Z$ to $X$.
 By Lemma~\ref{lem-4.2},
 we have $H^2(Z, T_Z(-F_1-F_2-F-D))=H^2(Z, T_Z(-\log(F_1+F_2+F+D)))=0$.
 Let $D_Z=F_1+F_2+F+D+S_1+S_2+S_3$.
 Since the self-intersection number of the section is $-1$,
 we still have the vanishing $H^2(Z, T_Z(-\log \ D_Z))=0$.
 We blow up eight times at the intersection points between
 three sections ($S_1$, $S_2$ and $S_3$) and two nodal fibers,
 and at the intersection points between
 two sections ($S_1$ and $S_2$) and one $I_2$-fiber.
 Denote this surface by $Z'$.
 Choose the exceptional curve in the total transform of $F_2$
 which intersects the proper transform of $S_1$, and choose two
 exceptional curves in the total transforms of $F_1$ and $F_2$
 which intersect the proper transform of $S_2$.
 Let $D_{Z'}$ be the reduced scheme of $F_1+F_2+F+D+S_1+S_2+S_3+$these
 three exceptional divisors.
 Then, by Lemma~\ref{lem-4.2}, Proposition~\ref{pro-2.6},
 and the self-intersection number, $-1$, of each exceptional divisor,
 we have $H^2(Z', T_{Z'}(-\log \ D_{Z'}))=0$.
 Finally, by using the same argument finite times through blowing up,
 we have the vanishing $H^2(\tilde Z, T_{\tilde Z}(-\log \ D_{\tilde Z}))=0$,
 where $D_{\tilde Z}$ are the five disjoint linear chains of $\mP^1$ which are
 the exceptional divisors from the contraction from $\tilde Z$ to $X$.
 Hence there is a $\mQ$-Gorenstein smoothing for $X$ by Theorem~\ref{thm-2.1}.
\end{proof}

\b

\section{Properties of $X_t$}
\label{sec-5}

 We showed in Section 4 that the projective surface $X$ has
 a $\mQ$-Gorenstein smoothing.
 We denote a general fiber of the $\mQ$-Gorenstein smoothing by $X_t$.
 In this section, we prove that $X_t$ is a simply connected, minimal,
 surface of general type with $p_g=0$ and $K_{X_t}^2=2$
 by using a standard argument.

 We first prove that $X_t$ satisfies $p_g=0$ and $K^2=2$:
 Since $\tilde Z$ is a nonsingular rational surface and $X$ has only
 rational singularities,
 $X$ is a projective surface with $H^1(X, \cO_X)=H^2(X, \cO_X)=0$.
 Then the upper semi-continuity implies
 $H^2(X_t,\cO_{X_t})=0$, so that the Serre duality implies $p_g(X_t)=0$
 (Equivalently, it also follows from the fact that
 $p_g(X_t)=\frac{b_2^+(X_t)-1}{2}=\frac{b_2^+(\tilde Z)-1}{2}=0$).
 And $K^2_X=2$ can be computed by using the explicit description of $f^*K_X$
 (refer to Equation~(\ref{eqn-2}) below).
 Then we have $K^2_{X_t}=2$ by the property of the $\mQ$-Gorenstein smoothing.

 Next, let us show the minimality of $X_t$: As we noticed in Section 3,
 the surface $\t Z$ contains the following two chains of $\mP^1$ including
 the proper transforms of three sections.
 We denote them by the following dual graphs
 $$\underset{G_1}{\overset{-2}{\circ}}-\underset{G_2}{\overset{-10}{\circ}}
 -\underset{G_3}{\overset{-2}{\circ}}
 -\underset{G_4}{\overset{-2}{\circ}}
   -\underset{G_5}{\overset{-2}{\circ}} -\underset{G_6}{\overset{-2}{\circ}}
   -\underset{G_7}{\overset{-2}{\circ}}-\underset{G_8}{\overset{-3}{\circ}},
   \ \ \
   \underset{J_1}{\overset{-5}{\circ}}-\underset{J_2}{\overset{-2}{\circ}}$$
 and we also denote four special fibers by the following dual graphs

\[\begin{array}{ccccccc}
 & \overset{\tilde A, -2}{\circ} &  & & \overset{J_2, -2}{\circ} - &
 \overset{J, -2}{\circ} &\\
  & \parallel &                     & &                        &    \vert  & \\
 \underset{E'_3, -1}{\circ}-& \underset{\tilde B, -4}{\circ} & -\underset{E''_3,
 -1}{\circ}, \ \ \ &
\underset{G_3, -2}{\circ} - & \underset{G_4, -2}{\circ} - &
\underset{G_5, -2}{\circ} & - \underset{G_6, -2}{\circ} -
\underset{G_7, -2}{\circ} ,
 \end{array}\]
\[\begin{array}{ccc}
 & {\overset{E'''_1, -1}{\circ}}  \\
 & \vert  & \\
 \underset{E'_1, -1}{\circ}- &\underset{I_1, -7}{\circ}&  -\underset{I_2,
 -2}{\circ}-\underset{I_3, -2}{\circ}-\underset{I_4,
-2}{\circ}-\underset{E''_1, -1}{\circ}\\
& \parallel & \\
& {\underset{E_1, -1}{\circ}}, &
\end{array}\]
\[\begin{array}{ccc}
 & {\overset{E'''_2, -1}{\circ}}&  \\
 & \vert  & \\
 \underset{E'_2, -1}{\circ}- \underset{H_1, -2}{\circ}- &
 \underset{H_2, -7}{\circ}& -\underset{H_3,
 -2}{\circ} -\underset{H_4, -2}{\circ}-\underset{H_5, -3}{\circ}
 -\underset{E''_2, -1}{\circ}-\underset{G_1, -2}{\circ}\\
 & \parallel & \\
 & {\underset{E_2, -1}{\circ}} . &
\end{array}\]
 Note that the second one indicates the $\tilde E_6$-singular fiber
 and $J$ denotes the rational $-2$-curve which is not contracted from
 $\tilde Z$ to $X$.
 The numbers indicate the self intersection numbers of curves.
 Let $f: \t Z\to X$ and let $h: \t Z\to Y=E(1)$.
 Then we have
\m
\[\begin{array}{ll}
  K_{\t Z} \equiv &
     f^*K_X - (\frac{7}{15}G_1+\frac{14}{15}G_2+\frac{13}{15}G_3
     +\frac{12}{15}G_4+\frac{11}{15}G_5+ \frac{10}{15}G_6+\frac{9}{15}G_7
     +\frac{8}{15}G_8+\frac{4}{9}H_1 \\
  &+\frac{8}{9}H_2+\frac{7}{9}H_3+\frac{6}{9}H_4+\frac{5}{9}H_5
   +\frac{4}{5}I_1+\frac{3}{5}I_2+\frac{2}{5}I_3+\frac{1}{5}I_4+\frac{1}{2}\tilde
   B+\frac{2}{3}J_1 +\frac{1}{3}J_2), \\
 \\
 K_{\t Z} \equiv &
 h^*K_Y+E_1+E'_1+4E''_1+E'''_1+I_2+2I_3+3I_4+E_2+2E'_2+8E''_2+E'''_2 \\
 &+H_1+H_3+2H_4+3H_5+4G_1+E'_3+E''_3
\end{array} \]
 On the other hand,
\[\begin{array}{ll}
 h^*K_Y \equiv &
 -\frac{1}{2}(2E_1+E'_1+E''_1+E'''_1+I_1+I_2+I_3+I_4)\\
 &-\frac{1}{2}(2E_2+E'_2+2E''_2+E'''_2+H_1+H_2+H_3+H_4+H_5+G_1)
\end{array}\]
\m
 Hence, combining these relations, we get

\begin{eqnarray}
\label{eqn-2}
\begin{array}{ll}
 f^*K_X  \equiv &
  \frac{119}{30}G_1+\frac{14}{15}G_2+\frac{13}{15}G_3+\frac{12}{15}G_4
  +\frac{11}{15}G_5+\frac{10}{15}G_6+\frac{9}{15}G_7+\frac{8}{15}G_8
  +\frac{17}{18}H_1 \\
  &+\frac{7}{18}H_2+\frac{23}{18}H_3+\frac{39}{18}H_4+\frac{55}{18}H_5
  +\frac{3}{10}I_1+\frac{11}{10}I_2+\frac{19}{10}I_3+\frac{27}{10}I_4+
  \frac{1}{2}\tilde B\\
  &+\frac{2}{3}J_1+\frac{1}{3}J_2+\frac{1}{2}E'_1+\frac{7}{2}E''_1+
  \frac{1}{2}E'''_1 +\frac{3}{2}E'_2+7E''_2+\frac{1}{2}E'''_2+E'_3+E''_3.
\end{array}
\end{eqnarray}

\m

\n Since all coefficients are positive in the expression of
 $f^*K_X$, the $\mQ$-divisor $f^*K_X$ is nef if $f^*K_X\cdot
 E_i'\ge 0$, $f^*K_X\cdot E_i''\ge 0$ for $i=1,2,3$,
 and $f^*K_X\cdot E_i'''\ge 0$ for $i=1,2$ - We have
 $f^*K_X\cdot E_1'=\frac{7}{15}$, $f^*K_X\cdot
 E_1''=\frac{2}{15}$, $f^*K_X\cdot
 E_1'''=\frac{5}{15}$, $f^*K_X\cdot E_2'=\frac{1}{9}$,
 $f^*K_X\cdot E_2''=\frac{1}{45}$,
 $f^*K_X\cdot E_2'''=\frac{19}{45}$, $f^*K_X\cdot E_3'=\frac{1}{6}$,
 and $f^*K_X\cdot E_3''=\frac{13}{30}$.
 Note that other divisors are contracted under the map $f$.
 The nefness of $f^*K_X$ implies the nefness of $K_X$.
 Since all coefficients are positive in the expression of
 $f^*K_X$, we get the vanishing $h^0(-K_X)=0$.
 Hence, by the upper semi-continuity property, i.e.
 the vanishing $h^0(-K_X)=0$ implies that $h^0(-K_{X_t})=0$,
 we conclude that $X_t$ is not a rational surface: If $X_t$ is a
 rational surface with $h^0(-K_{X_t})=0$, then $\chi(2K_{X_t})\le 0$.
 But $\chi(2K_{X_t})=\chi(\cO_{X_t})+K_{X_t}^2=3$,
 which is a contradiction.
 Since $K_{X_t}^2=2$, $X_t$ is a surface of general type
 by the classification theory of surfaces.
 Let $\pi: \cX \to \Delta$ be a $\mQ$-Gorenstein smoothing of $X$.
 Since the $\mQ$-Cartier divisor $K_{\cX/\Delta}$ is $\pi$-big over
 $\Delta$ and $\pi$-nef at the point 0, the nefness of $K_{X_t}$
 is also obtained by shrinking $\Delta$ if it is necessary \cite{Nak}.
 Therefore we have

\begin{proposition}
\label{prop-5.1}
 $X_t$ is a minimal surface of general type with
 $p_g=0$ and $K_{X_t}^2=2$.
\end{proposition}

 Finally, applying the standard arguments about Milnor fibers (\cite{LW}, \S5),
 we prove Theorem~\ref{thm-main}: Note that $X$ has five
 singularities, say $p_1, \ldots, p_5$.
 For each $1 \leq i \leq 5$, let $(V_i, p_i)$ be a small disjoint
 neighborhood for the singularity $p_i$ in $X$.
 By Proposition~\ref{pro-ksb},
 there exists a closed embedding of $(V_i, p_i)$ in
 $(\mC^n, 0)$ satisfying that a local $\mQ$-Gorenstein smoothing $V$
 of $V_i$ is a closed embedding in $(\mC^n\times\Delta, 0)$ and
 the map $\pi: V\to\Delta$ is induced by
 the projection on the second factor $\mC^n\times\Delta\to\Delta$.
 Let $B_r=\{ z\in \mC^n \ | \ |z|<r\}$ and let $S_r=\partial B_r$.
 The sphere $S_r$ is called a Milnor sphere for $V_i$ if for every $0< r'\le r$
 the sphere $S_{r'}$ intersects $V_i$ transversally.
 Let $S_r$ be a Milnor sphere for $V_i$ then we can assume that
 $S_r \times \Delta$ intersects $V_t$
 transversally for all $t \in \Delta$ by shrinking $\Delta$ if it is
 necessary. In this set-up $F_i:=V_t\cap (B_r \times \Delta)$
 is called the Minor fiber of a $\mQ$-Gorenstein smoothing.
 Then there are small disjoint neighborhoods $(V_i, p_i) \ (1 \leq i \leq 5)$
 for the five singularities $p_i$ in $X$ such that
 $\hat{X}:= (X-\cup_{i=1}^5 V_i)\cup(\cup_{i=1}^5 F_i) \ \
 \text{is diffeomorphic to} \ \
 X_t =(X-\cup_{i=1}^5 F_i)\cup(\cup_{i=1}^5 F_i),$
 where the pasting is made by choosing an orientation preserving
 diffeomorphism $\partial V_i\to \partial F_i$ for each $1 \leq i \leq 5$.
 Note that, since these $F_i$ are rational balls (\cite{Man01}, \S1),
 the manifold $\hat{X}=(X-\cup_{i=1}^5 V_i)\cup(\cup_{i=1}^5 F_i)$
 is diffeomorphic to the rational blow-down $4$-manifold $\t Z_{15, 9, 5, 3, 2}$
 constructed in Theorem~\ref{thm-3.1}.
 Hence the simple connectivity of $X_t$ follows from the fact
 that $\t Z_{15, 9, 5, 3, 2}$ is simply connected.

\b

\section{More examples}
\label{sec-6}

 In this section we construct more examples of simply connected, minimal,
 complex surfaces of general type with $p_g=0$ and $K^2 =2$ by using
 different configurations coming from different elliptic pencils in $\mP^2$.
 Since all the proofs are basically the same as the case of the main example
 constructed in Section 3, we only explain how to construct more examples.

\m

\begin{figure}[hbtb]
 \begin{center}
 \setlength{\unitlength}{1mm}
 \includegraphics[height=3.5cm]{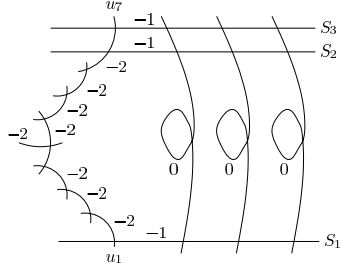}
 \end{center}
 \caption{A rational surface $Y'$}
 \label{fig-Y2}
\end{figure}

\n {\bf Construction.} We first consider an elliptic fibration
 on $E(1)$ which has one $\tilde E_7$-singular fiber ($III^*$ in Kodaira's
 table of singular elliptic fibers in~\cite{BHPV}, p.201), three nodal fibers
 and three disjoint sections. Such an elliptic fibration can constructed
 explicitly as follows (\cite{ss} for details):
 Let $C_1$ be a union of two lines $L_1$ and $L_2$ with the latter
 of multiplicity two, and let $C_2$ be an immersed $2$-sphere with
 one positive transverse double point.
 Assume that $L_2$ intersects $C_2$ at a single point, say $p$,
 with triple tangency between two curves, and $L_1$ (also passing
 through $p$) intersects $C_2$ at two other smooth points, say $q$
 and $r$.
 Consider a cubic pencil in $\mP^2$ induced by $C_1$ and $C_2$,
 i.e. $\lambda C_1+\mu C_2$, for $[\lambda:\mu]\in \mP^1$.
 Now, blowing up seven times successively at the point $p$ and
 blow up once at the points $q$ and $r$ respectively,
 we get a desired elliptic fibration $E(1)=\mP^2\sharp 9\overline{\mP}^2$
 over $\mP^1$. We denote this elliptic fibration by $g: Y'=E(1) \to \mP^1$.
 Note that there is an $\tilde E_7$-singular fiber coming from the
 total transforms of $L_1$ and $L_2$, and three nodal fibers on $Y'$.
 Furthermore, the elliptic fibration $Y'$ admits three disjoint sections,
 say $S_1, S_2, S_3$, where $S_1$ is the $7^{th}$-exceptional
 curve and $S_2, S_3$ are the $8^{th}$- and the $9^{th}$-exceptional
 curves in the forming of $Y'$ (Figure~\ref{fig-Y2}).
 Notice that the section $S_1$ meets an ending $2$-sphere, say $u_1$, of the
 $\tilde E_7$-singular fiber, and two sections $S_2$ and $S_3$ meet
 the other ending $2$-sphere, say $u_7$, of the $\tilde E_7$-singular fiber.

\begin{figure}[hbtb]
 \begin{center}
 \setlength{\unitlength}{1mm}
 \includegraphics[height=4cm]{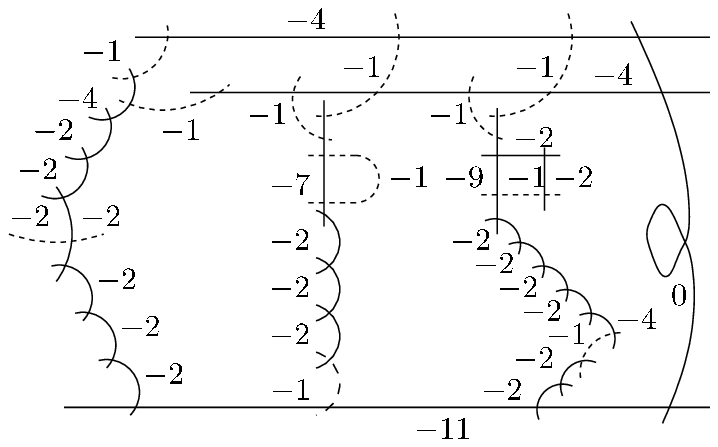}
 \end{center}
 \caption{A rational surface $\t Z'$}
 \label{fig-tZ2}
\end{figure}

 Let $Z':=Y'\sharp 2{\overline \mP}^2$ be the surface obtained
 by blowing up at two singular points of two nodal fibers on $Y'$.
 Then
 there are two fibers such that each consists of two $\mP^1$s,
 say $E_i$ and $F_i$, satisfying $E_i^2=-1$, $F_i^2=-4$ and
 $E_i\cdot F_i=2$ for $i=1,2.$
 Note that each $E_i$ is an exceptional curve and
 $F_i$ is the proper transform of a nodal fiber.
 First, we blow up six times at the intersection points
 between two sections $S_2, S_3$ and $F_1, F_2, u_7$.
 It makes the self-intersection number of the proper transforms
 of $S_2$, $S_3$ and $u_7$ to be $-4$, and the self-intersection
 number of the proper transforms of $F_1$ and $F_2$ are $-6$.
 We also blow up four times successively at the intersection point
 between $F_1$ and the section $S_1$ which makes a chain of $\mP^1$,
 ${\overset{-7}{\circ}}-{\overset{-2}{\circ}}-{\overset{-2}{\circ}}
 -{\overset{-2}{\circ}}$, and blow up six times
 successively at the intersection point between $F_2$ and $S_1$
 which makes a chain of $\mP^1$,
 ${\overset{-7}{\circ}}-{\overset{-2}{\circ}}-{\overset{-2}{\circ}}
 -{\overset{-2}{\circ}}-{\overset{-2}{\circ}}-{\overset{-2}{\circ}}$.
 We further blow up twice successively at one intersection point
 between $-7$-curve in the total transform of $F_2$ and the exceptional
 curve intersecting $-7$-curve twice, and finally blow up twice successively
 at the intersection point between the ending $-2$-curve and the exceptional
 curve in the total transform of $F_2$, so that it makes a chain of $\mP^1$,
 ${\overset{-2}{\circ}}-{\overset{-2}{\circ}}-{\overset{-9}{\circ}}
 -{\overset{-2}{\circ}}-{\overset{-2}{\circ}}-{\overset{-2}{\circ}}
 -{\overset{-2}{\circ}}-{\overset{-4}{\circ}}$.
 Hence we have a rational surface $\t Z':= Y' \sharp 22{\overline \mP}^2$
 which contains five disjoint linear chains of $\mP^1$:
 ${\overset{-2}{\circ}}-{\overset{-2}{\circ}}-{\overset{-11}{\circ}}
 -{\overset{-2}{\circ}}-{\overset{-2}{\circ}}-{\overset{-2}{\circ}}
 -{\overset{-2}{\circ}}-{\overset{-2}{\circ}}-{\overset{-2}{\circ}}
 -{\overset{-4}{\circ}}$
 (which contains the proper transforms of $S_1$, the six rational
 $-2$-curves in $\tilde E_7$-singular fiber and the proper transform of $u_7$),
 ${\overset{-2}{\circ}}-{\overset{-2}{\circ}}-{\overset{-9}{\circ}}
 -{\overset{-2}{\circ}}-{\overset{-2}{\circ}}-{\overset{-2}{\circ}}
 -{\overset{-2}{\circ}}-{\overset{-4}{\circ}}$
 (which is a part of the total transform of $F_2$),
 ${\overset{-7}{\circ}}-{\overset{-2}{\circ}}-{\overset{-2}{\circ}}
 -{\overset{-2}{\circ}}$ (which is a part of the total transform of
 $F_1$), and two $-4$-curves which are the proper transforms of
 $S_2$ and $S_3$ (Figure~\ref{fig-tZ2}).

 Next, we contract these five disjoint chains of $\mP^1$ from $\t Z'$.
 This creates a projective surface $X$ with five singularities of class $T$.
 The existence of a $\mQ$-Gorenstein smoothing for $X$ is obtained
 by the same proof as in Section 4.
 Let us denote a general fiber of the $\mQ$-Gorenstein smoothing by $X_t$.
 Then, by the same argument in Section 5, we see that
 $X_t$ is a minimal surface of general type with $p_g=0$ and $K^2=2$.
 Finally, since the rational blow-down $4$-manifold ${\t Z'}_{25,19,5,2,2}$
 of the surface $\t Z'$ is simply connected, which can be proved in a similar
 way as in the proof of Theorem~\ref{thm-3.1},
 we conclude that $X_t$ is a simply connected, minimal, complex surface
 of general type with $p_g=0$ and $K^2 =2$.

\m

\n {\bf Remark.} J. Keum suggested that one can also construct
 another example of simply connected surface of general type with
 $p_g=0$ and $K^2=2$ by using an elliptic fibration on $E(1)$ which
 has one $\tilde D_6$-singular fiber, one $I_2$-singular fiber and
 two nodal fibers, and which also admits four disjoint sections~\cite{Ke}.
 Of course, the remaining proofs are the same as in Section 4 and 5.
 We do not know whether all these constructions above provide
 the same deformation equivalent type of surfaces with $p_g=0$ and $K^2=2$.
 It is a very intriguing problem to determine whether the examples
 constructed above are diffeomorphic (or deformation equivalent)
 to the surface constructed in Section 3.

\b

\section{Appendix: Simply connected surfaces with $p_g=0$ and $K^2=1$}
\label{sec-appendix}

 As mentioned in the Introduction, the only previously known simply
 connected complex surface of general type with $p_g=0$ was Barlow surface,
 which has $K^2 =1$. Dolgachev-Werner~\cite{DW}
 constructed a minimal surface of general type with $p_g=0$ and
 $K^2=1$ by using a special quintic surface with one elliptic singular
 point, which was founded by Craighero-Gattazzo~\cite{CG}.
 Their surface has a trivial first homology group,
 but the simple connectivity is still unknown.
 More examples and systematic study of numerical Godeaux surfaces, i.e.
 the minimal surfaces of general type with $p_g=0$ and $K^2=1$,
 are given in~\cite{CP, reid06}.
 In this appendix we construct simply connected, minimal, complex
 surfaces of general type with $p_g=0$ and $K^2 =1$ using a rational
 blow-down surgery and a $\mQ$-Gorenstein smoothing theory.
 Since all proofs are basically the same as the case of $K^2=2$,
 we only explain how to construct such examples.

\m

\begin{figure}[hbtb]
 \begin{center}
 \setlength{\unitlength}{1mm}
 \includegraphics[height=3.5cm]{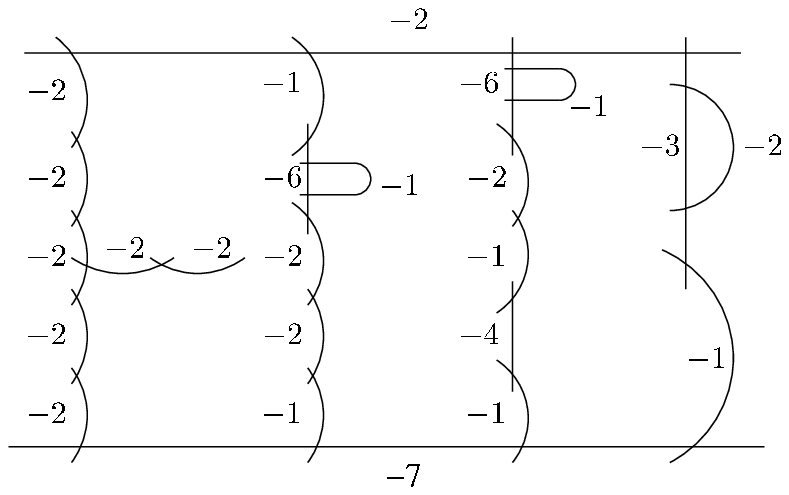}
 \end{center}
 \caption{A rational surface $\t Z$}
 \label{fig-tZ3}
\end{figure}

\n {\bf Construction A1.}
 We begin with the rational elliptic surface
 $Y=\mP^2\sharp 9\overline{\mP}^2$ with special fibers used in Section 3.
 We use the same notations for $Z, F_1, F_2, E_1, E_2, S_1, S_2,
 \tilde A, \tilde B$ as in the main construction in Section 3.
 First, we blow up at the intersection points between
 $F_1$ and the section $S_1$, so that the self-intersection number
 of the proper transform of $S_1$ is $-2$.
 We also blow up at the intersection points between $F_1, F_2, \tilde B$
 and the section $S_2$. It makes the self-intersection number of the proper
 transform of $S_2$ to be $-4$, and the self-intersection number of the
 proper transform of $\tilde B$ to be $-3$.
 Next,
 we blow up twice successively at the intersection point between the proper
 transform of $S_2$ and the exceptional curve in the total transform of $F_1$,
 so that it makes a chain of $\mP^1$,
 ${\overset{-6}{\circ}}-{\overset{-2}{\circ}}-{\overset{-2}{\circ}}$,
 lying in the total transform of $F_1$.
 And we blow up at the intersection point
 between the proper transform of $F_2$ and the exceptional curve
 which intersects the proper transform of the section $S_2$.
 That makes a chain of $\mP^1$,
 ${\overset{-6}{\circ}}-{\overset{-1}{\circ}}-{\overset{-2}{\circ}}$,
 lying in the total transform of $F_2$.
 Then we blow up again at the intersection point
 between the exceptional $-1$-curve and the rational $-2$-curve, and
 blow up again at the intersection point between the rational
 $-3$-curve,  which is the proper transform of the rational $-2$-curve,
 and the proper transform of the section $S_2$.
 It makes a chain of $\mP^1$,
 ${\overset{-6}{\circ}}-{\overset{-2}{\circ}}-{\overset{-1}{\circ}}-
 {\overset{-4}{\circ}}-{\overset{-1}{\circ}}$.
 Note that this process makes the self-intersection number of the proper
 transform of $S_2$ to be $-7$.
 Finally, we have a rational surface $\t Z:= Y \sharp 11{\overline \mP}^2$
 which contains four disjoint linear chains of $\mP^1$:
 ${\overset{-7}{\circ}}-{\overset{-2}{\circ}}-{\overset{-2}{\circ}}
 -{\overset{-2}{\circ}}$ (which contains the proper transform of $S_2$
 and a part of $\tilde E_6$-fiber),
 ${\overset{-6}{\circ}}-{\overset{-2}{\circ}}-{\overset{-2}{\circ}}$,
 ${\overset{-2}{\circ}}-{\overset{-6}{\circ}}
 -{\overset{-2}{\circ}}-{\overset{-3}{\circ}}$
 (which contains the proper transforms of $S_1$ and of $\tilde B$),
 and ${\overset{-4}{\circ}}$ (Figure~\ref{fig-tZ3}).

 Next, by contracting these four disjoint chains of $\mP^1$ from $\t Z$,
 we get a projective surface $X$ with four singularities of class $T$.
 Then, by the same argument in Section 4 and 5 above,
 we finally see that the surface $X$ has a $\mQ$-Gorenstein smoothing
 and a general fiber $X_t$ of the $\mQ$-Gorenstein smoothing for $X$
 is a simply connected, minimal, complex surface of general type with
 $p_g=0$ and $K^2=1$.

\m

\begin{figure}[hbtb]
 \begin{center}
 \setlength{\unitlength}{1mm}
 \includegraphics[height=3.5cm]{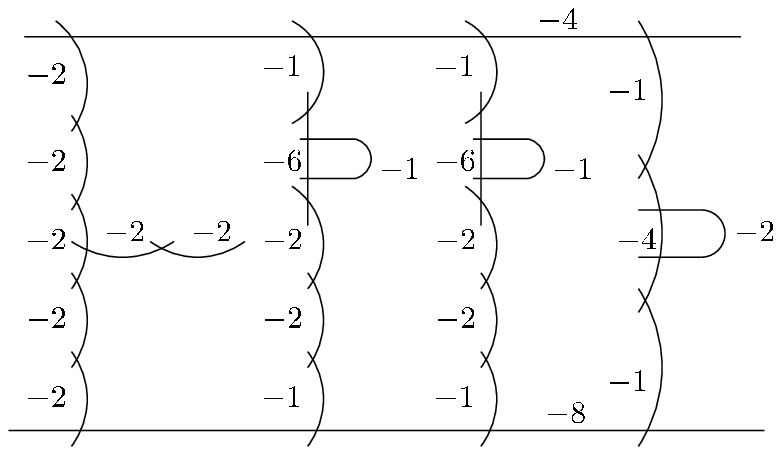}
 \end{center}
 \caption{A rational surface $\t Z'$}
 \label{fig-tZ4}
\end{figure}

\n {\bf Construction A2.}
 Similar to Construction A1 above, we use the same elliptic fibration on $Y$
 and we use the same notations for $Z, F_1, F_2, E_1, E_2, S_1, S_2,
 \tilde A, \tilde B$ as in the main construction in Section 3.
 First, we blow up at the intersection points between
 $F_1, F_2, \tilde B$ and the sections $S_1, S_2$ respectively,
 so that the self-intersection numbers
 of the proper transforms of $S_1$ and $S_2$ are $-4$.
 And we blow up twice successively at each intersection point between
 the proper transform of $S_2$ and the exceptional curves in the total
 transforms of $F_1, F_2$, so that it makes two disjoint linear chains
 of $\mP^1$,
 ${\overset{-6}{\circ}}-{\overset{-2}{\circ}}-{\overset{-2}{\circ}}$
 and
 ${\overset{-6}{\circ}}-{\overset{-2}{\circ}}-{\overset{-2}{\circ}}$,
 lying in the total transforms of $F_1$ and $F_2$ respectively.
 Note that this process makes the self-intersection number of the proper
 transform of $S_2$ to be $-8$.
 Then, using four consecutive $2$-spheres in $\tilde E_6$-fiber which
 is not connected to the section $S_1$ together with the proper
 transform of $S_2$, we also get a linear chain of
 $\mP^1$, ${\overset{-8}{\circ}}-{\overset{-2}{\circ}}-{\overset{-2}{\circ}}
 -{\overset{-2}{\circ}}-{\overset{-2}{\circ}}$.
 Hence we finally have a rational surface $\t Z':= Y \sharp 12{\overline \mP}^2$
 which contains five disjoint linear chains of $\mP^1$:
 ${\overset{-8}{\circ}}-{\overset{-2}{\circ}}-{\overset{-2}{\circ}}
 -{\overset{-2}{\circ}}-{\overset{-2}{\circ}}$
 (which contains the proper transform of $S_2$ and a part of $\tilde E_6$-fiber),
 ${\overset{-6}{\circ}}-{\overset{-2}{\circ}}-{\overset{-2}{\circ}}$,
 ${\overset{-6}{\circ}}-{\overset{-2}{\circ}}-{\overset{-2}{\circ}}$,
 ${\overset{-4}{\circ}}$ (which is the proper transform  of $S_1$)
 and ${\overset{-4}{\circ}}$ (which is the proper transform of $\tilde B$)
 (Figure~\ref{fig-tZ4}).

 Then, as the same way in the Construction A1 above,
 we get a simply connected, minimal, complex surface of general type
 with $p_g=0$ and $K^2=1$.

\m

\n {\bf Open Problem.} Determine whether the examples
 constructed in the Appendix above are diffeomorphic
 (or deformation equivalent) to the Barlow surface.

\b
\b

\end{document}